\documentclass[12pt]{article}
\usepackage{latexsym,geometry,cite}
\newtheorem{theorem}{Theorem}

\newtheorem{corollary}[theorem]{Corollary}

\usepackage{lineno}
\usepackage{setspace}

\begin{document}

\onehalfspace

\title{Extremal Values of the Chromatic Number\\ for a Given Degree Sequence}

\author{St\'{e}phane Bessy$^1$ and Dieter Rautenbach$^2$}

\date{}

\maketitle

\begin{center}
$^1$ 
Laboratoire d'Informatique, de Robotique et de Micro\'{e}lectronique de Montpellier (LIRMM), Montpellier, France, \texttt{stephane.bessy@lirmm.fr}\\[3mm]
$^2$ Institute of Optimization and Operations Research, Ulm University, Ulm, Germany, \texttt{dieter.rautenbach@uni-ulm.de}
\end{center}

\begin{abstract}
For a degree sequence $d:d_1\geq \cdots \geq d_n$,
we consider the smallest chromatic number $\chi_{\min}(d)$ and the largest chromatic number $\chi_{\max}(d)$
among all graphs with degree sequence $d$.
We show that if $d_n\geq 1$, then $\chi_{\min}(d)\leq \max\left\{ 3,d_1-\frac{n+1}{4d_1}+4\right\}$,
and, if $\sqrt{n+\frac{1}{4}}-\frac{1}{2}>d_1\geq d_n\geq 1$, then $\chi_{\max}(d)=\max\limits_{i\in [n]}\min\left\{ i,d_i+1\right\}$.
For a given degree sequence $d$ with bounded entries,
we show that 
$\chi_{\min}(d)$,
$\chi_{\max}(d)$,
and also
the smallest independence number $\alpha_{\min}(d)$ among all graphs with degree sequence $d$,
can be determined in polynomial time.
\end{abstract}

{\small 

\begin{tabular}{lp{13cm}}
{\bf Keywords:} & Degree sequence; chromatic number; independence number\\
{\bf MSC 2010:} & 05C07; 05C15; 05C69
\end{tabular}
}

\pagebreak

\section{Introduction}

We consider finite, simple, and undirected graphs.
The {\it degree sequence} of a graph $G$ with vertex set $\{ v_1,\ldots,v_n\}$ is the sequence $d_G(v_1),\ldots,d_G(v_n)$ of its vertex degrees.
A sequence $d_1,\ldots,d_n$ of integers is a {\it degree sequence} if it is the degree sequence of some graph.
Repetitions within the degree sequence can be indicated by suitable exponents; 
the degree sequence of the star $K_{1,r}$ of order $r+1$, for instance, is $r,1^r$.
For a given sequence $d$, let ${\cal G}(d)$ be the set of all graphs $G$ whose degree sequence is $d$;
called the {\it realizations} of $d$. 
For an integer $n$, let $[n]$ be the set of the positive integers at most $n$.

In the present paper we consider
$$\chi_{\min}(d)=\min\left\{\chi(G):G\in {\cal G}(d)\right\}
\,\,\,\,\,\,\,\,\,\,\,\,\,\mbox{and}\,\,\,\,\,\,\,\,\,\,\,\,\,
\chi_{\max}(d)=\max\left\{\chi(G):G\in {\cal G}(d)\right\}.$$
Punnim \cite{pu} determined $\chi_{\min}(d)$ and $\chi_{\max}(d)$ for regular degree sequences $d=r^n$ in almost all cases.
The parameter $\chi_{\max}(d)$ was also considered by Dvo\v{r}\'{a}k and Mohar \cite{dvmo},
who established degree sequence versions of the Hadwiger Conjecture and even the Haj\'{o}s Conjecture, see also \cite{roso}.

We contribute some bounds, exact values, and algorithmic results.
Further discussion of related research will be given throughout the rest of the paper.

\section{Some bounds and exact values}

For a sequence $d$ of non-negative integers $d_1\geq \cdots\geq d_n$, 
let $H(d)$ be the sequence $$d_2-1,\ldots,d_{d_1+1}-1,d_{d_1+2},\ldots,d_n.$$
Havel \cite{hav} and Hakimi \cite{hak} showed that $d$ is a degree sequence if and only if $H(d)$ is a degree sequence.
In fact, they observed that if $d$ is a degree sequence,  
then there is a realization $G$ of $d$ in which the neighbours of a vertex of degree $d_1$ have degrees $d_2,\ldots,d_{d_1+1}$.
Iteratively applying this observation to a given degree sequence yields a realization 
that tends to contain a large complete subgraph on the vertices of large degrees,
that is, such a realization may be expected to have high chromatic number.

In order to obtain a realization with hopefully small chromatic number,
one can apply Havel and Hakimi's observation to the complement.
More precisely, for a degree sequence $d$ as above, the sequence $\bar{d}$ defined as 
$$n-1-d_n\geq\cdots\geq n-1-d_1$$ 
is also a degree sequence;
in fact, the graphs in ${\cal G}\left(\bar{d}\right)$ are exactly the complements $\bar{G}$ of the graphs $G$ in ${\cal G}(d)$.
Furthermore, by the above observation of Havel and Hakimi,
$\bar{d}$ has a realization in which the neighbors of a vertex of the largest degree $n-1-d_n$ have degrees $n-1-d_{n-1},\ldots,n-1-d_{d_n+1}$.
Equivalently, as already observed by Kleitman and Wang \cite{klwa} in a more general form,
$d$ has a realization in which the neighbors of a vertex of the smallest degree $d_n$ have degrees $d_1,\ldots,d_{d_n}$.
In summary, we obtain that $d$ is a degree sequence if and only if the sequence $\bar{H}(d)$ defined as 
\begin{eqnarray}\label{ehbar}
d_1-1,\ldots,d_{d_n}-1,d_{d_n+1},\ldots,d_{n-1}
\end{eqnarray}
is a degree sequence.
Iteratively applying this observation to a given degree sequence yields a realization 
that tends to avoid dense subgraphs on the vertices of large degrees,
that is, such a realization may be expected to have small chromatic number.

As an example consider the degree sequence $d:r^{r+1},1^{r(r+1)}$ for some positive integer $r$.
Havel and Hakimi's original observation yields the realization $K_{r+1}\cup {r\choose 2}K_2$, 
whose chromatic number is $r+1$, which equals $\chi_{\max}(d)$,
while the above complementary version yields the realization $(r+1)K_{1,r}$, 
whose chromatic number is $2$, which equals $\chi_{\min}(d)$.

\medskip

\noindent For a sequence $d$ of integers $d_1,\ldots,d_n$, 
let $n$ be the {\it length} of $d$, 
let $\min(d)=\min\{ d_1,\ldots,d_n\}$, 
and let $\max(d)=\max\{ d_1,\ldots,d_n\}$.
Furthermore, let $\bar{H}^0(d)=d$, $\bar{H}^1(d)=\bar{H}(d)$, and 
$\bar{H}^i(d)=\bar{H}\left(\bar{H}^{i-1}(d)\right)$ for an integer $i$ at least $2$.
Note that iteratively applying the reductions $d\mapsto H(d)$ or $d\mapsto\bar{H}(d)$ 
always requires reordering the constructed sequences in a non-increasing way.

\begin{theorem}\label{theorem3}
If $d$ is a degree sequence of length $n$, then
$$\chi_{\min}(d)\leq \max\left\{\min\left(\bar{H}^{n-i}(d)\right):i\in [n]\right\}+1.$$
\end{theorem}
{\it Proof:}
Iteratively applying the complementary version of Havel and Hakimi's observation to the degree sequence $d$
yields a realization $G$ of $d$ with vertex set $\{ v_1,\ldots,v_n\}$
such that, for $i$ from $n$ down to $1$,
the vertex $v_i$ has degree $\min\left(\bar{H}^{n-i}(d)\right)$ in the graph $G[\{ v_1,\ldots,v_i\}]$.
Greedily coloring the vertices of $G$ in the order $v_1,\ldots,v_n$ 
yields a coloring that uses at most $\max\left\{\min\left(\bar{H}^{n-i}(d)\right):i\in [n]\right\}+1$ colors.
$\Box$

\medskip

\noindent Note that for the degree sequence $d:r^{r+1},1^{r(r+1)}$ of length $n=(r+1)^2$ considered as an example above,
we obtain $\max\left\{\min\left(\bar{H}^{n-i}\left(r^{r+1},1^{r(r+1)}\right)\right):i\in [n]\right\}+1=2$,
that is, for this degree sequence $d$, Theorem \ref{theorem3} reproduces the correct value of $\chi_{\min}(d)$.

Unfortunately, Theorem \ref{theorem3} is not very explicit. 
As a more explicit consequence, we quantify how small degrees may reduce the effect of large degrees on $\chi_{\min}(d)$.

\begin{corollary}\label{corollary1}
If $d$ is a degree sequence $d_1\geq \ldots\geq d_n$, 
and $k$ and $\ell$ are positive integers such that $d_k\geq k+\ell$ and $d_{n-\ell+1}\leq k$, then
$$\chi_{\min}(d)\leq \max\left\{ d_1-\frac{1}{k}\left(1+\sum\limits_{i=n-\ell+1}^n d_i\right)+1,d_{k+1},k\right\}+1.$$
\end{corollary}
{\it Proof:} 
We consider the first $\ell$ applications of the reduction $d\mapsto\bar{H}(d)$.
Since $d_k\geq k+\ell$ and $d_{n-\ell+1}\leq k$, 
we obtain that, for $i\in [\ell]$,
the degree sequence $\bar{H}^i(d)$ arises from $\bar{H}^{i-1}(d)$ by removing the degree $d_{n-i+1}$,
and reducing the $d_{n-i+1}$ largest degrees by $1$.
For $i\in \{ 0,\ldots,\ell\}$, let $\Delta_i=\max\left(\bar{H}^i(d)\right)$, 
and let $n_i$ be the number of entries of $\bar{H}^i(d)$ that are equal to $\Delta_i$.
Suppose, for a contradiction, that $\Delta_{\ell}>\max\left\{ d_1-\frac{D+1}{k}+1,d_{k+1}\right\}$,
where $D=\sum\limits_{i=n-\ell+1}^n d_i$.
Note that each of the $\ell+1$ degree sequences $d,\bar{H}(d),\ldots,\bar{H}^{\ell}(d)$ 
contains at most $k$ entries that are strictly larger than $d_{k+1}$.
So, for $i\in [\ell]$, we have
\begin{itemize}
\item $(\Delta_i,n_i)=(\Delta_{i-1},n_{i-1}-d_{n-i+1})$ if $d_{n-i+1}<n_{i-1}$, and
\item $\Delta_i=\Delta_{i-1}-1$ and $n_i\leq k-(d_{n-i+1}-n_{i-1})=n_{i-1}-d_{n-i+1}+k$ if $d_{n-i+1}\geq n_{i-1}$.
\end{itemize}
Note that $(k\Delta_{i-1}+n_{i-1})-(k\Delta_i+n_i)\geq d_{n-i+1}$ in both cases.
Summation over $i\in [\ell]$ yields $(k\Delta_0+n_0)-(k\Delta_{\ell}+n_{\ell})\geq D$.
Since $\Delta_0=d_1$, $n_0\leq k$, and $n_{\ell}\geq 1$, this implies 
$\Delta_{\ell}\leq d_1-\frac{D+1}{k}+1$,
which is a contradiction.
Hence, $\Delta_{\ell}\leq\max\left\{ d_1-\frac{D+1}{k}+1,d_{k+1}\right\}$, and 
any realization $H$ of the degree sequence $\bar{H}^{\ell}(d)$ 
can be colored using at most $\max\left\{ d_1-\frac{D+1}{k}+1,d_{k+1}\right\}+1$ many colors.
Adding $\ell$ further vertices of degrees $d_{n-\ell+1},\ldots,d_n$ one by one to $H$,
and connecting them to suitable vertices according to the previous reductions, yields a realization $G$ of $d$.
Since the added vertices all have degree at most $k$,
the coloring of $H$ can be extended greedily to a coloring of $G$ 
using at most $\max\left\{ d_1-\frac{D+1}{k}+1,d_{k+1},k\right\}+1$ different colors in total.
$\Box$

\medskip

\noindent For a given degree sequence $d$ not satisfying any further restriction, 
one can only bound $\chi_{\min}(d)$ from above by $\max(d)+1$.
In fact, $d$ might be $\max(d)^{\max(d)+1},0^{n-\max(d)-1}$, 
whose only realization contains a clique of size $\max(d)+1$.

Our next two results improve this trivial estimate 
for graphs without isolated vertices.

\begin{theorem}\label{theorem5}
If $d$ is a degree sequence of length $n$ with $\max(d)\geq\sqrt{\frac{n\delta}{4}}$ and $\min(d)\geq \delta$ for some positive integer $\delta$,
then $\chi_{\min}(d)\leq \max(d)-\frac{n\delta}{4\max(d)}+\delta+3$.
\end{theorem}
{\it Proof:} Our first goal is to show that we may assume that $d$ has a realization with a very large independent set.
Therefore, among all realizations $G$ of the degree sequence $d$ and all (not necessarily optimal) colorings $f$ of $G$, 
we choose $G$ and $f$ with color classes $V_1,\ldots,V_k$, where $V_i$ contains $n_i$ vertices for $i\in [k]$,
in such a way that 
\begin{itemize}
\item $(n_1,\ldots,n_k)$ is lexicographically maximal, and
\item subject to this first condition, the number of edges between $V_{k-1}$ and $V_k$ is minimum.
\end{itemize}
Note that $k$ may actually be larger than $\chi(G)$,
and that $n_1$ is necessarily equal to the independence number $\alpha(G)$ of $G$.

Let $\Delta=\max(d)$.
If $k\leq \Delta-\frac{n\delta}{4\Delta}+\delta+3$, then $\chi_{\min}(d)\leq \chi(G)\leq k$ implies the desired bound.
Hence, we may assume that $k>\Delta-\frac{n\delta}{4\Delta}+\delta+3$.
Since $\Delta\geq \sqrt{\frac{n\delta}{4}}$ and $\delta\geq 1$, we have $k\geq 5$. 
By the choice of the coloring $f$, there is an edge, say $uv$, between the smallest two color classes $V_{k-1}$ and $V_k$.
If $G\setminus (V_{k-1}\cup V_k\cup N_G(u)\cup N_G(v))$ contains an edge $xy$,
then removing from $G$ the two edges $uv$ and $xy$,
and adding the two edges $ux$ and $vy$, 
yields another realization $G'$ of $d$.
Note that $f$ is still a coloring of $G'$.
This implies that there is a coloring $f'$ of $G'$ such that
either the non-increasing vector of the sizes of the color classes
is lexicographically larger than the one of $f$,
or there are fewer edges between the two smallest color classes.
Since both cases imply a contradiction to the choice of $G$ and $f$,
we obtain that $V(G)\setminus (V_{k-1}\cup V_k\cup N_G(u)\cup N_G(v))$ is an independent set, which implies
$\alpha(G)\geq n-(n_{k-1}+n_k)-2\Delta$.
Since $V_{k-1}$ and $V_k$ are the smallest two color classes, and $n_2+\cdots+n_k=n-\alpha(G)$,
we obtain $n_{k-1}+n_k\leq \frac{2}{k-1}(n-\alpha(G))$.
This implies $\alpha(G)\geq n-\frac{2}{k-1}(n-\alpha(G))-2\Delta$,
and, using $k\geq 5$, we obtain $\alpha(G)\geq n-\frac{k-1}{k-3}\cdot 2\Delta\geq n-4\Delta$.

Altogether, we may assume that $d$ has a realization $G$ 
with an independent set $I=\{ u_1,\ldots,u_{\alpha}\}$ of order at least $n-4\Delta$.
By the above-mentioned observations of 
Havel \cite{hav}, Hakimi \cite{hak}, Rao \cite{ra1}, and Kleitman and Wang \cite{klwa},
we may further assume that, for every $i\in [\alpha]$,
the vertex $u_i$ is adjacent to $d_G(u_i)$ vertices in $V(G)\setminus I$ 
of the largest degrees in the induced subgraph $G-\{ u_1,\ldots,u_{i-1}\}$ of $G$.
Arguing as in the proof of Corollary \ref{corollary1}, we obtain
$\Big((n-\alpha)\Delta+(n-\alpha)\Big)-\Big((n-\alpha)\Delta(G-I)+1\Big)\geq d_G(u_1)+\cdots+d_G(u_{\alpha})\geq \alpha\delta$,
where $\Delta(G-I)$ denotes the maximum degree of $G-I$.
This implies 
$\Delta(G-I)
\leq \Delta-\frac{\alpha\delta+1}{n-\alpha}+1
\leq \Delta-\frac{(n-4\Delta)\delta+1}{4\Delta}+1
= \Delta-\frac{n\delta+1}{4\Delta}+\delta+1$.
Therefore, we can color $G$ using at most $\Delta-\frac{n\delta+1}{4\Delta}+\delta+2$ colors on the vertices in $V(G)\setminus I$,
and one additional color on the vertices in $I$,
which implies $\chi_{\min}(d)\leq \chi(G)\leq \Delta-\frac{n\delta+1}{4\Delta}+\delta+3$.
$\Box$

\medskip

\noindent For positive integers $r$, $s$, and $\delta$ such that $r+1$ is a multiple of $\delta$,
let $d$ be the degree sequence $(r+s)^{r+1},\delta^{s(r+1)/\delta}$.
Since the sum of the largest $r+1$ degrees equals exactly 
$2{r+1\choose 2}+\delta s(r+1)/\delta$,
every realization $G$ of $d$ contains a clique on the $r+1$ vertices of largest degrees,
and an independent set on the remaining vertices.
Note that $\chi(G)\in \{ r+1,r+2\}$,
which, for $r\gg s\gg\delta$, is roughly 
$\max(d)-\frac{n\min(d)}{\max(d)}$,
that is, up to the constants, the bound in Theorem \ref{theorem5} is best possible.
In fact, by imposing a stronger lower bound on $\max(d)$
or by increasing the additive constant,
the factor $4$ within the term $\frac{n\delta+1}{4\Delta}$
can easily be reduced to slightly more than $2$.

Our next result gives a best possible bound on $\chi_{\min}(d)$ for degree sequences of small degrees.

\begin{theorem}\label{theorem1}
If $n, d_1,\ldots,d_n$ are integers such that 
$\sqrt{\frac{n-1}{2}}\geq d_1\geq \cdots \geq d_n\geq 1$ and $d_1+\cdots+d_n$ is even, then $\chi_{\min}(d)\leq 3$.
(In particular, $d_1,\ldots,d_n$ is a degree sequence.)
\end{theorem}
{\it Proof:} 
There is a partition of $[n]$ into two sets $X$ and $Y$ with 
$\Big| |X|-|Y| \Big|\leq 1$ and 
$0\leq s\leq d_1\leq \sqrt{\frac{n-1}{2}}$,
where $s=\sum\limits_{i\in X}d_i-\sum\limits_{i\in Y}d_i$;
in fact, as long as there are two equal entries $d_i$ and $d_j$ in the sequence $d_1,\ldots,d_n$,
we assign $i$ to $X$ and $j$ to $Y$, and remove $d_i$ and $d_j$ from the sequence,
and once all remaining entries are distinct, say $d_{i_1}>\cdots >d_{i_k}$, 
we assign $i_1,i_3,\ldots$ to $X$ and $i_2,i_4,\ldots$ to $Y$.
Let $x=|X|$ and $y=|Y|$.
Note that $x,y\geq \frac{n-1}{2}$; in particular, $s\leq x$.
Reducing $s$ distinct entries of the sequence $(d_i)_{i\in X}$ by $1$, and reordering yields a sequence $a_1\geq \cdots \geq a_x$.
Reordering the sequence $(d_i)_{i\in Y}$ yields $b_1\geq \cdots \geq b_y$.

By construction, $\sum\limits_{i\in [x]}a_i=\sum\limits_{i\in [y]}b_i$,
$\max\{ a_1,b_1\}\leq \sqrt{\frac{n-1}{2}}$, and $b_y\geq 1$.

Let $k\in [x]$. 
If $k\leq \sqrt{\frac{n-1}{2}}$, then $a_1\leq \sqrt{\frac{n-1}{2}}$ and $b_n\geq 1$ imply
$$\sum\limits_{i\in [k]}a_i\leq ka_1\leq \frac{n-1}{2}\leq y\leq \sum\limits_{i\in [y]}\min\{ k, b_i\}.$$
If $k>\sqrt{\frac{n-1}{2}}$, then $b_1\leq \sqrt{\frac{n-1}{2}}$ implies
$$\sum\limits_{i\in [k]}a_i\leq \sum\limits_{i\in [x]}a_i=\sum\limits_{i\in [y]}b_i=\sum\limits_{i\in [y]}\min\{ k, b_i\}.$$
By the Gale-Ryser Theorem \cite{ga,ry}, 
there is a bipartite graph $H$ with partite sets $X$ and $Y$ with $|X|=x$ and $|Y|=y$
such that the vertices in $X$ have degrees $a_1,\ldots,a_x$
and the vertices in $Y$ have degrees $b_1,\ldots,b_y$.
Since $s$ has the same parity as $\sum\limits_{i\in X}d_i+\sum\limits_{i\in Y}d_i=d_1+\cdots+d_n$, it is an even integer,
and adding to $H$ a matching of size $s/2$ incident to those vertices in $X$ 
corresponding to the entries of $(d_i)_{i\in X}$ that were previously reduced by $1$,
results in a graph $G$ with degree sequence $d_1,\ldots,d_n$.
Clearly, $\chi(G)\leq 3$, and the upper bound on $\chi_{\min}(d)$ follows.
$\Box$

\medskip

\noindent The conclusion of Theorem \ref{theorem1} is best possible, 
because there might not be a subset $X$ of $[n]$ with $\sum\limits_{i\in X}d_i=\sum\limits_{i\in [n]\setminus X}d_i$,
which is a necessary condition for the existence of a bipartite realization.  
The complexity of deciding the existence of a bipartite realization for a given degree sequence is unknown.

Note that together, Theorem \ref{theorem5} and Theorem \ref{theorem1} imply
$$\chi_{\min}(d)\leq \max\left\{ 3,\max(d)-\frac{n+1}{4\max(d)}+4\right\}$$
for every degree sequence $d$ with $\min(d)\geq 1$.

Theorem \ref{theorem1} has the following variant where the essential assumption is that $\max(d)-\min(d)$ is small.
Note that this next result also covers regular degree sequences of sufficient length.

\begin{theorem}\label{theorem1b}
If $n, d_1,\ldots,d_n$ are integers and $\epsilon >0$ is such that 
$\frac{n-1}{2}\epsilon\ge d_1\ge \cdots \ge d_n\ge 1$,
$d_1-d_n\leq \sqrt{\frac{n-1}{2}}(1-\epsilon)$, and
$d_1+\cdots+d_n$ is even, then $\chi_{\min}(d)\leq 3$.
\end{theorem}
{\it Proof:} We may assume that $d_1>\sqrt{\frac{n-1}{2}}$; otherwise Theorem \ref{theorem1} implies the result.
Furthermore, we have $\epsilon\leq 1$.
Exactly as in
the proof of Theorem \ref{theorem1}, we obtain the existence of a
partition of $[n]$ into two sets $X$ and $Y$ with $\Big| |X|-|Y|
\Big|\leq 1$ and $0\leq s\leq d_1\leq \frac{n-1}{2}\epsilon$, where
$s=\sum\limits_{i\in X}d_i-\sum\limits_{i\in Y}d_i$.  Setting $x=|X|$
and $y=|Y|$, we obtain, as above, that $x,y\geq \frac{n-1}{2}$, $s\leq
x$, and $s$ is even.  Let $a_1\geq \cdots \geq a_x$ and $b_1\geq
\cdots \geq b_y$ be as in the proof of Theorem \ref{theorem1}.  By
construction, $\sum\limits_{i\in [x]}a_i=\sum\limits_{i\in [y]}b_i$,
$\max\{ a_1,b_1\}\leq d_1$, and $b_y\geq d_n$.

Notice that as $d_1>\sqrt{\frac{n-1}{2}}$, we have
$$\frac{d_n}{d_1}\ge \frac{d_1-\sqrt{\frac{n-1}{2}}(1-\epsilon)}{d_1}\ge 1-(1-\epsilon)=\epsilon.$$
Let $k\in [x]$.
If $k\leq d_n$, then 
$$\sum\limits_{i\in [k]}a_i\leq kd_1\leq k\frac{n-1}{2}\leq ky\leq \sum\limits_{i\in [y]}\min\{ k, b_i\}.$$
If $d_n<k<d_1$, then 
$$\sum\limits_{i\in [k]}a_i\leq kd_1\leq d_1^2 \leq \frac{n-1}{2} \epsilon d_1\le
\frac{n-1}{2}d_n\leq yd_n\leq
\sum\limits_{i\in [y]}\min\{ k, b_i\}.$$
And, if $k\geq d_1$, then
$$\sum\limits_{i\in [k]}a_i\leq \sum\limits_{i\in [x]}a_i=\sum\limits_{i\in [y]}b_i=\sum\limits_{i\in [y]}\min\{ k, b_i\}.$$
At this point, the proof can be completed exactly as the proof of Theorem \ref{theorem1}. $\Box$

\medskip

\noindent For a graph $G$ with degree sequence $d_1\geq \cdots \geq d_n$, 
Welsh and Powell \cite{wepo} observed
\begin{eqnarray}\label{ewp}
\chi(G)\leq \max\limits_{i\in [n]}\min\left\{ i,d_i+1\right\},
\end{eqnarray}
which is an immediate consequence of applying the natural greedy coloring algorithm to the vertices of $G$ in an order of non-increasing degrees.
If $d_1\geq\cdots\geq d_n$ is a degree sequence
such that $d_p-d_{p+1}\geq p-2$ for $p=\max\limits_{i\in [n]}\min\left\{ i,d_i+1\right\}$,
then Havel and Hakimi's observation explained above implies the existence of a realization $G$ of $d$ 
for which the vertices of degrees $d_1,\ldots,d_p$ form a clique.
This implies $p\leq \chi(G)\leq \chi_{\max}(d)\leq p$, that is, 
$\chi_{\max}(d)=\max\limits_{i\in [n]}\min\left\{ i,d_i+1\right\}$ for such degree sequences.

Our next result shows that the Welsh-Powell bound (\ref{ewp}) also gives the correct value of $\chi_{\max}(d)$ 
for degree sequences $d$ of small degrees.

\begin{theorem}\label{theorem2}
If $n, d_1,\ldots,d_n$ are integers such that 
$\sqrt{n+\frac{1}{4}}-\frac{1}{2}>d_1\geq \cdots \geq d_n\geq 1$ and $d_1+\cdots+d_n$ is even, 
then $\chi_{\max}(d)=\max\limits_{i\in [n]}\min\left\{ i,d_i+1\right\}$.
\end{theorem}
{\it Proof:} Let $p=\max\limits_{i\in [n]}\min\left\{ i,d_i+1\right\}$. 
Note that $p\leq d_p+1\leq d_1+1$.

By the Welsh-Powell bound (\ref{ewp}), every graph $G$ with degree sequence $d_1,\ldots,d_n$ satisfies $\chi(G)\leq p$,
which implies $\chi_{\max}(d)\leq p$.
In order to establish equality, 
we show the existence of a realization that contains a clique of size $p$.

Let $k\in [n]$.
We obtain 
$\sum\limits_{i\in [k]}d_i\leq k d_1$
and 
$k(k-1)+\sum\limits_{i\in [n]\setminus [k]}\min\{ k,d_i\}\geq k(k-1)+n-k$.
Therefore, $\sum\limits_{i\in [k]}d_i$ is at most $k(k-1)+\sum\limits_{i\in [n]\setminus [k]}\min\{ k,d_i\}$
if $k d_1\leq k(k-1)+n-k$, which is equivalent to $k(d_1+2-k)\leq n$.
Since 
$\sqrt{n+\frac{1}{4}}-\frac{1}{2}>d_1\geq 1$
implies 
$n\geq 3$
and 
$k(d_1+2-k)\leq \left(\frac{d_1+2}{2}\right)^2\leq n$,
the Erd\H{o}s-Gallai Theorem \cite{erga} implies the existence of a graph with degree sequence $d_1,\ldots,d_n$.
Among all such graphs with vertex set $\{ v_1,\ldots,v_n\}$, 
where $v_i$ has degree $d_i$ for $i\in [n]$,
we choose $G$ such that the number $m(G[\{ v_1,\ldots,v_p\}])$
of edges of the subgraph of $G$ induced by $\{ v_1,\ldots,v_p\}$ is as large as possible.

Suppose, for a contradiction, that $G[\{ v_1,\ldots,v_p\}]$ is not a clique, that is,
$v_i$ and $v_j$ are not adjacent in $G$ for distinct $i$ and $j$ in $[p]$.
By the choice of $p$, we have $d_i,d_j\geq p-1$, 
which implies that $v_i$ and $v_j$ both have at least one neighbor in $R=\{ v_{p+1},\ldots,v_n\}$.

First, we assume that $v_i$ and $v_j$ both have the same unique neighbor $v_r$ in $R$,
that is, $\{ v_r\}=N_G(v_i)\cap R=N_G(v_j)\cap R$.
Since there are at most $1+d_1^2$ vertices at distance at most $2$ from $v_r$, including, in particular, $v_i$ and $v_j$,
and $n-(p-2)-(1+d_1^2)\geq n-d_1^2-d_1>0$,
there is a vertex $v_s$ in $R$ with a neighbor $v_t$ such that $v_s$ and $v_t$ are both not adjacent to $v_r$.
Now, 
removing from $G$ the edges $v_iv_r$, $v_jv_r$, and $v_sv_t$, and 
adding the edges $v_iv_j$, $v_rv_s$, and $v_rv_t$ 
yields a realization $G'$ of $d_1,\ldots,d_n$ with $m(G'[\{ v_1,\ldots,v_p\}])>m(G[\{ v_1,\ldots,v_p\}])$,
which contradicts the choice of $G$.

Now, we may assume that $v_i$ is adjacent to some vertex $v_r$ in $R$, and that $v_j$ is adjacent to a different vertex $v_s$ in $R$.
If $v_r$ is not adjacent to $v_s$, 
then removing from $G$ the edges $v_iv_r$ and $v_jv_s$, and 
adding the edges $v_iv_j$ and $v_rv_s$ 
yields a realization $G'$ of $d_1,\ldots,d_n$ with $m(G'[\{ v_1,\ldots,v_p\}])>m(G[\{ v_1,\ldots,v_p\}])$,
which contradicts the choice of $G$.
Hence, we may assume that $v_r$ and $v_s$ are adjacent.
Since there are at most $1+d_1^2$ vertices at distance at most $2$ from $v_r$, including, in particular, $v_i$, $v_s$, and $v_j$,
and $n-(p-2)-(1+d_1^2)\geq n-d_1^2-d_1>0$,
there is a vertex $v_p$ in $R$ with a neighbor $v_q$ such that 
$v_p$ is not adjacent to $v_s$,
and 
$v_q$ is not adjacent to $v_r$.
Note that $v_q$ may be $v_j$, in which case, $v_j$ has distance $2$ from $v_r$.
Now, 
removing from $G$ the edges $v_iv_r$, $v_jv_s$, and $v_pv_q$, and 
adding the edges $v_iv_j$, $v_sv_p$, and $v_rv_q$ 
yields a realization $G'$ of $d_1,\ldots,d_n$ with $m(G'[\{ v_1,\ldots,v_p\}])>m(G[\{ v_1,\ldots,v_p\}])$,
which contradicts the choice of $G$.

Altogether, we obtain that $G$ contains a clique of order $p$, which completes the proof. $\Box$

\section{Algorithmic aspects}

One way to establish that $\chi_{\max}(d)$ is large is to show the existence of a realization of $d$ that contains a large clique.
Dvo\v{r}\'{a}k and Mohar \cite{dvmo} proved the best possible statement 
that for every degree sequence $d$, some realization of $d$ has a clique of size at least $5/6(\chi_{\max}(d)-3/5)$.
Since Rao \cite{ra1,ra2} efficiently characterized the largest clique size $\omega_{\max}(d)$ of any realization of a given degree sequence $d$,
and, trivially, $\chi_{\max}(d)\geq \omega_{\max}(d)$,
we immediately obtain that $\chi_{\max}(d)$ can be approximated in polynomial time for a given $d$ within an asymptotic factor of $6/5$. 

Our next two results show that $\chi_{\max}(d)$ and $\chi_{\min}(d)$ can both be determined in polynomial time
for given degree sequences with bounded entries.

\begin{corollary}\label{corollary2}
Let $\Delta$ be a fixed positive integer.

For a given degree sequence $d$ with $\max(d)\leq \Delta$,
one can determine $\chi_{\max}(d)$ in polynomial time.
\end{corollary}
{\it Proof:} Let $d$ have length $n$.
Clearly, we may assume $\min(d)\geq 1$.
If $\sqrt{n-2}\geq \Delta$, then Theorem \ref{theorem2} implies that $\chi_{\max}(d)$ coincides with the Welsh-Powell bound (\ref{ewp}).
If $\sqrt{n-2}<\Delta$, then, as $\Delta$ is fixed, there are only constantly many realizations of $d$,
which can all be generated and optimally colored by brute force in constant time. $\Box$

\begin{theorem}\label{theorem4}
Let $k$ and $p$ be fixed positive integers.

For a given degree sequence $d$ with at most $p$ distinct entries,
one can decide in polynomial time whether $\chi_{\min}(d)\leq k$.
\end{theorem}
{\it Proof:} Let $d:d_1^{n_1},\ldots,d_p^{n_p}$ and $n=n_1+\cdots+n_p$.
There are 
$\prod\limits_{i=1}^p{n_i+k-1\choose k-1}\leq \left(\frac{n}{p}+k\right)^{kp}$
distinct matrices $(n_i^j)_{(i,j)\in [p]\times [k]}$ with non-negative integral entries $n_i^j$
such that $\sum\limits_{j=1}^kn_i^j=n_i$ for $i\in [p]$.
It is easy to see that $\chi_{\min}(d)\leq k$ if and only if there is such a matrix $(n_i^j)_{(i,j)\in [p]\times [k]}$ 
for which the complete $k$-partite graph whose $j$th partite set $V_j$ has order $\sum\limits_{i=1}^p n_i^j$ for $j\in [k]$,
has a factor $G$ such that $V_j$ contains exactly $n_i^j$ vertices of degree $d_i$ in $G$ for every $i\in [p]$ and $j\in [k]$.
Since the existence of such a factor can be decided in polynomial time using matching methods,
and, for fixed $k$ and $p$, there are only polynomially many different suitable matrices, 
the desired statement follows. $\Box$

\medskip

\noindent It seems plausible to wonder whether $\chi_{\max}(d)$ is linked to $\alpha_{\min}(d)$,
the minimum independence number of a realization of $d$.
While $\alpha_{\max}(d)=\omega_{\max}\left(\bar{d}\right)$
can be determined efficiently using the results of Rao \cite{ra1,ra2},
Bauer, Hakimi, Kahl, and Schmeichel \cite{bahakasc} conjectured that it is computationally hard to determine $\alpha_{\min}(d)$ 
for a given degree sequence $d$.

Our next goal is to show that also $\alpha_{\min}(d)$ can be determined in polynomial time for given degree sequences $d$ with bounded entries.

For a degree sequence $d_1,\ldots,d_n$, let $\alpha_{CW}(d)=\sum\limits_{i=1}^n\frac{1}{d_i+1}$.
Caro \cite{ca} and Wei \cite{we} proved that $\alpha(G)\geq \alpha_{CW}(d)$
for every graph $G$ with degree sequence $d$.
For a connected graph $G$ with degree sequence $d$, 
Harant and Rautenbach \cite{hara} showed 
$\alpha(G)\geq k\geq \sum\limits_{u\in V(G)}\frac{1}{d_G(u)-f(u)+1}$,
where $k$ is an integer, and, for every vertex $u$ of $G$, $f(u)$ is a non-negative integer at most $d_G(u)$ such that
$\sum\limits_{u\in V(G)}f(u)\geq 2(k-1)$.
This improved an earlier result of Harant and Schiermeyer \cite{hasc}.

If $\alpha_{CW}(d)\geq 2$, then $k\geq \alpha_{CW}(d)$ implies $2(k-1)\geq k\geq \alpha_{CW}(d)$, and, hence,
\begin{eqnarray*}
\alpha(G) & \geq & \sum\limits_{u\in V(G)}\frac{1}{d_G(u)-f(u)+1}\\
&=& \alpha_{CW}(d)+\sum\limits_{u\in V(G)}\left(\frac{1}{d_G(u)-f(u)+1}-\frac{1}{d_G(u)+1}\right)\\
&\geq & \alpha_{CW}(d)+\frac{1}{(\max(d)+1)^2}\sum\limits_{u\in V(G)}f(u)\\
&\geq & \left(1+\frac{1}{(\max(d)+1)^2}\right)\alpha_{CW}(d).
\end{eqnarray*}

\begin{theorem}\label{theorem6}
Let $\Delta$ be a fixed positive integer.

For a given degree sequence $d$ with $\max(d)\leq \Delta$,
every component of every realization $G$ of $d$ with $\alpha(G)=\alpha_{\min}(d)$ 
has order at most $\left((\Delta+1)^3+1\right)\left(\left(\frac{\Delta+2}{2}\right)^2+{\Delta+1\choose 2}\right)$.
In particular, one can determine $\alpha_{\min}(d)$ in polynomial time.
\end{theorem}
{\it Proof:}
Let $d$ be a degree sequence with $\max(d)\leq \Delta$.
Let $G$ be a realization of $d$ with $\alpha(G)=\alpha_{\min}(d)$.
Suppose, for a contradiction, that some component $K$ of $G$ has order $n(K)$ more than the stated value.
Let $R$ be a set of $\left(\frac{\Delta+2}{2}\right)^2$ vertices of $K$.
For $i\in [\Delta]$, let $V_i$ be the set of vertices of degree $i$ in $V(K)\setminus R$, and let $n_i=|V_i|$.
Let $p_i=\left\lfloor\frac{n_i}{i+1}\right\rfloor$, and let $S_i$ arise by removing $p_i(i+1)$ vertices from $V_i$ for each $i\in [\Delta]$.
Note that $|S|\leq \sum\limits_{i=1}^\Delta i={\Delta+1\choose 2}$, where $S=S_1\cup\cdots\cup S_{\Delta}$,
that is, $R\cup S$ is a set of at least $\left(\frac{\Delta+2}{2}\right)^2$
and at most $\left(\frac{\Delta+2}{2}\right)^2+{\Delta+1\choose 2}$
many vertices of $K$.
Let $d'$ be the sequence of the degrees of the vertices in $R\cup S$,
and let $d''$ be the sequence of the degrees of the vertices in $V(K)\setminus (R\cup S)$.
Note that $\alpha_{CW}(d'')\geq \frac{(n(K)-|R\cup S|)}{\Delta+1}$.
Hence, the lower bound on $n(K)$ implies
$\left(1+\frac{1}{(\Delta+1)^2}\right)\alpha_{CW}(d'')
=\frac{1}{(\Delta+1)^2}\alpha_{CW}(d'')+\alpha_{CW}(d'')
>|R\cup S|+\alpha_{CW}(d'')$.
As observed in the proof of Theorem \ref{theorem2},
the Erd\H{o}s-Gallai Theorem implies that the sequence $d'$,
which is a sequence of positive integers at most $\Delta$ that is of length at least $\left(\frac{\Delta+2}{2}\right)^2$,
is a degree sequence.
Let $K'_0$ be a realization of $d'$.
By construction, the graph $K'=K'_0\cup\bigcup\limits_{i=1}^{\Delta}p_iK_{i+1}$ has exactly the same degree sequence as $K$.
By the result of Harant and Rautenbach mentioned above,
\begin{eqnarray*}
\alpha(K') & = & \alpha(K'_0)+\sum\limits_{i=1}^{\Delta}p_i\alpha(K_{i+1})\\
& = & \alpha(K'_0)+\alpha_{CW}(d'')\\
& \leq & |R\cup S|+\alpha_{CW}(d'')\\
& < & \left(1+\frac{1}{(\Delta+1)^2}\right)\alpha_{CW}(d'')\\
& < & \left(1+\frac{1}{(\Delta+1)^2}\right)\alpha_{CW}(d)\\
& \leq & \alpha(K).
\end{eqnarray*}
Therefore, replacing $K$ by $K'$ within $G$ 
yields a realization $G'$ of $d$ with $\alpha(G')<\alpha(G)$,
contradicting the choice of $G$.
This completes the proof of the first part of the statement.

Since, as $\Delta$ is fixed, there are only finitely many graphs of maximum degree at most $\Delta$
and order at most $\left((\Delta+1)^3+1\right)\left(\left(\frac{\Delta+2}{2}\right)^2+{\Delta+1\choose 2}\right)$.
Listing, for each of these graphs, the degree sequence and the independence number,
it is a routine matter to determine $\alpha_{\min}(d)$ 
for a given degree sequence $d$ with $\max(d)\leq \Delta$
by dynamic programming in polynomial time.
$\Box$

\end{document}